\documentclass[12pt]{article}

\usepackage{amssymb}
\usepackage{amsfonts}
\usepackage{latexsym}
\usepackage{psfig}

\newcommand{\bi}{\bibitem}
\newcommand{\ci}{\cite}

\newtheorem{theo}{Theorem}[section]
\newtheorem{lemma}{Lemma}[section]

\begin{document}


\def\a{\alpha}           \def\cA{{\cal A}}     \def\bA{{\bf A}}
\def\b{\beta}            \def\cB{{\cal B}}     \def\bB{{\bf B}}
\def\g{\gamma}           \def\cC{{\cal C}}     \def\bC{{\bf C}}
\def\G{\Gamma}           \def\cD{{\cal D}}     \def\bD{{\bf D}}
\def\d{\delta}           \def\cE{{\cal E}}     \def\bE{{\bf E}}
\def\D{\Delta}           \def\cF{{\cal F}}     \def\bF{{\bf F}}
\def\ve{\varepsilon}     \def\cG{{\cal G}}     \def\bG{{\bf G}}
\def\z{\zeta}            \def\cH{{\cal H}}     \def\bH{{\bf H}}
\def\e{\eta}             \def\cI{{\cal I}}     \def\bI{{\bf I}}
\def\vt{\vartheta}       \def\cJ{{\cal J}}     \def\bJ{{\bf J}}
\def\vT{\Theta}          \def\cK{{\cal K}}     \def\bK{{\bf K}}
\def\k{\kappa}           \def\cL{{\cal L}}     \def\bL{{\bf L}}
\def\l{\lambda}          \def\cM{{\cal M}}     \def\bM{{\bf M}}
\def\L{\Lambda}          \def\cN{{\cal N}}     \def\bN{{\bf N}}
\def\m{\mu}              \def\cO{{\cal O}}     \def\bO{{\bf O}}
\def\n{\nu}              \def\cP{{\cal P}}     \def\bP{{\bf P}}
\def\r{\rho}             \def\cQ{{\cal Q}}     \def\bQ{{\bf Q}}
\def\s{\sigma}           \def\cR{{\cal R}}     \def\bR{{\bf R}}
\def\S{\Sigma}           \def\cS{{\cal S}}     \def\bS{{\bf S}}
\def\t{\tau}             \def\cT{{\cal T}}     \def\bT{{\bf T}}
\def\f{\phi}             \def\cU{{\cal U}}     \def\bU{{\bf U}}
\def\F{\Phi}             \def\cV{{\cal V}}     \def\bV{{\bf V}}
\def\vp{\varphi}         \def\cW{{\cal W}}     \def\bW{{\bf W}}
\def\c{\chi}             \def\cX{{\cal X}}     \def\bX{{\bf X}}
\def\p{\psi}             \def\cY{{\cal Y}}     \def\bY{{\bf Y}}
\def\P{\Psi}             \def\cZ{{\cal Z}}     \def\bZ{{\bf Z}}
\def\o{\omega}
\def\O{\Omega}
\def\x{\xi}
\def\X{\Xi}


\def\Z{{\mathbb Z}}
\def\R{{\mathbb R}}
\def\C{{\mathbb C}}
\def\T{{\mathbb T}}
\def\N{{\mathbb N}}
\def\S{{\mathbb S}}

\def\ma{\left(\begin{array}{cc}}
\def\am{\end{array}\right)}
\def\iint{\int\!\!\!\int}
\def\lt{\biggl}
\def\rt{\biggr}
\let\geq\geqslant
\let\leq\leqslant
\def\[{\begin{equation}}
\def\]{\end{equation}}
\def\wt{\widetilde}
\def\pa{\partial}
\def\sm{\setminus}
\def\es{\emptyset}
\def\no{\noindent}
\def\ol{\overline}
\def\iy{\infty}
\def\ev{\equiv}
\def\/{\over}
\def\ts{\times}
\def\os{\oplus}
\def\ss{\subset}
\def\h{\hat}
\def\Re{\mathop{\rm Re}\nolimits}
\def\Im{\mathop{\rm Im}\nolimits}
\def\supp{\mathop{\rm supp}\nolimits}
\def\sign{\mathop{\rm sign}\nolimits}
\def\Ran{\mathop{\rm Ran}\nolimits}
\def\Ker{\mathop{\rm Ker}\nolimits}
\def\Tr{\mathop{\rm Tr}\nolimits}
\def\const{\mathop{\rm const}\nolimits}
\def\Wr{\mathop{\rm Wr}\nolimits}
\def \BBox{\hspace{1mm}\vrule height6pt width5.5pt depth0pt \hspace{6pt}}


\def\Twelve{
\font\Tenmsa=msam10 scaled 1200
\font\Sevenmsa=msam7 scaled 1200
\font\Fivemsa=msam5 scaled 1200
\textfont\msbfam=\Tenmsb
\scriptfont\msbfam=\Sevenmsb
\scriptscriptfont\msbfam=\Fivemsb

\font\Teneufm=eufm10 scaled 1200
\font\Seveneufm=eufm7 scaled 1200
\font\Fiveeufm=eufm5 scaled 1200
\textfont\eufmfam=\Teneufm
\scriptfont\eufmfam=\Seveneufm
\scriptscriptfont\eufmfam=\Fiveeufm}

\def\Ten{
\textfont\msafam=\tenmsa
\scriptfont\msafam=\sevenmsa
\scriptscriptfont\msafam=\fivemsa

\textfont\msbfam=\tenmsb
\scriptfont\msbfam=\sevenmsb
\scriptscriptfont\msbfam=\fivemsb

\textfont\eufmfam=\teneufm
\scriptfont\eufmfam=\seveneufm
\scriptscriptfont\eufmfam=\fiveeufm}
\makeatletter
\def\eqalign#1{\null\vcenter{\def\\{\cr}\openup\jot\m@th
  \ialign{\strut$\displaystyle{##}$\hfil&$\displaystyle{{}##}$\hfil
      \crcr#1\crcr}}\,}
\makeatother

\title {Spectral estimates for periodic Jacobi matrices}

\author{Evgeni Korotyaev
\\ Institut f\"ur  Mathematik\\  Humboldt Universit\"at zu Berlin \\
e-mail:  ek@mathematik.hu-berlin.de
 \and Igor V. Krasovsky\\ 
 Institut f\"ur  Mathematik\\  Technische Universit\"at, Berlin \\e-mail:
ivk@math.tu-berlin.de }

\maketitle

\begin{abstract}
We obtain bounds for the spectrum and for the total width of 
the spectral gaps for Jacobi matrices on $\ell^2(\Z)$ of the form 
$(H\psi)_n=
a_{n-1}\psi_{n-1}+b_n\psi_n+a_n\psi_{n+1}$, where 
$a_n=a_{n+q}$ and $b_n=b_{n+q}$ are periodic sequences of real numbers.
The results are based on a study of the quasimomentum $k(z)$ corresponding 
to $H$. We consider $k(z)$ as a conformal mapping in the complex plane.
We obtain the trace identities which connect integrals of the 
Lyapunov exponent over the gaps with the normalised traces of powers of $H$. 

\end{abstract}

\vskip 0.25cm

\section {Introduction}
\setcounter{equation}{0}

The purpose of the present work is to study the spectrum of the 
$q$-periodic Jacobi matrix on $\ell^2(\Z)$
\[
(H\psi)_n=
a_{n-1}\psi_{n-1}+b_n\psi_n+a_n\psi_{n+1},\label{H}
\]
where
$b_n=b_{n+q}\in \R$, $a_n=a_{n+q}>0$, $n\in \Z$. The number $q$ is the 
{\it smallest} period.
(It is well known that any symmetric Jacobi matrix can be 
represented in the form with positive off-diagonal elements.)

Periodic Jacobi matrices were discussed in many works (see, e.g.,
\ci{Toda,Moerbeke,Perkolab}).

In what follows, we consider $q>1$ thus
excluding the trivial case 
of $q=1$ (the spectrum of $H$ with $q=1$ is the interval 
$[b_1-2a_1,b_1+2a_1]$). 
Let $r(H)$ denote half the
distance between the ends of the spectrum of $H$. We set 
\[
c=r(H),
\]
\[
A=(a_1a_2\cdots a_q)^{1/q}.
\]
We will also use the notation $a_j\equiv a$
for $a_1=a_2=\cdots=a_q=a$.
We remark that periodic 
discrete 1D Schr\"odinger operators are a particular case of (\ref{H}) with
$a_j\equiv 1$. 

By Gershgorin's theorem (see, e.g., \ci{Lancaster}), 
the spectrum of $H$ lies inside the interval 
$[\min_j(b_j-a_j-a_{j-1}),\max_j(b_j+a_j+a_{j-1})]$.
It is absolutely continuous and consists of 
{\it exactly} $q$ intervals
$$
\s_m=[\l^+_{m-1},\l^-_m],  \quad m=1,\dots,q
$$
separated by gaps
$$
\gamma_m=(\l^-_m,\l^+_m), \quad m=1,\dots,q-1.
$$
If a gap degenerates, i.e. $\g_m=\es $, then the corresponding segments
$\s_m$, $\s_{m+1}$ merge. In what follows, we will use $|\cdot |$ to denote 
the Lebesgue measure of sets.

It is known \ci{DS,Last} that the total width of bands and gaps, 
respectively, satisfy the inequalities\footnote{They were written in 
\ci{Last} for the Schr\"odinger case ($a_j\equiv 1$) but the arguments
are easy to generalise.}
\[
2c-\tilde b
\ge\sum_{m=1}^q|\s_m|>{4A^q\over M^{q-1}},\qquad
\sum_{m=1}^{q-1}|\gamma_m|\ge \tilde b,\label{eL}
\]
where $\tilde b=\max_j b_j-\min_j b_j$ and
$M=\max(\max_j(b_j+a_j+a_{j-1})-\min_jb_j,
\max_jb_j-\min_j(b_j-a_j-a_{j-1}))$. Analysis of the proof in \ci{Last} 
shows that $M$ in (\ref{eL}) can be replaced by $2c$.

In the present work we shall find further
estimates for the widths of the gaps 
and for the total width of the spectrum
in terms of the matrix elements of $H$.
Our results are obtained from properties of the quasimomentum
(associated with the periodic operator $H$)
considered as a conformal mapping of complex domains.
Originally, the method of conformal mappings was proposed in \ci{MO}
for the Hill operator ($Hy(x)=y''(x)+f(x)y(x)$, $f(x)=f(x+T)$). It was further 
developed in \ci{KK1,K1,K2,K3,K4} to obtain spectral estimates for 
the Hill operator and the Dirac operator.
As we shall see below,
the main peculiarity of the discrete case is that 
the spectrum is bounded. This makes $z=\arccos
(\l/\max|\l^\pm_n|)$
a more natural spectral variable than $\l$. Thus our estimates 
for the gaps will be formulated in terms of $z$.

As is known,
the real part of the quasimomentum is the integrated density
of states of the operator $H$ while the imaginary part is the
Lyapunov exponent (it describes the rate of growth in $n$ of the solutions 
$\psi_n(\l)$ of the equation $(H-\l)\psi=0$).
In this paper we shall obtain identities (trace formulas and Dirichlet
integrals) which relate various integrals of the quasimomentum
to traces of powers of $H$ (Lemmas \ref{Ltr} and \ref{Ldir}). 
These identities will serve as a starting point to derive our
estimates. They are also of separate interest.

We shall call a Jacobi matrix normalised if $-\l^+_0=\l^-_q>0$.
Obviously, this can be achieved for an arbitrary matrix by adding 
a certain constant.
Below we always assume that $H$ is normalised.

We shall need the matrix $L$ defined by
\begin{equation}
\eqalign{
L=\pmatrix{b_1 & a_1+ia_2\cr a_1-ia_2 & b_2}\quad \mathrm{if}\; q=2,\qquad 
L=\pmatrix{b_1 & a_1 &&&& ia_q\cr
a_1 & b_2 & a_2 &\cr
  & a_2 & b_3 & a_3 \cr
  &   &     & \ddots \cr
 &&& a_{q-2} & b_{q-1} &a_{q-1}\cr
-ia_q&&&&a_{q-1}&b_{q}},\\
\mathrm{if}\; q>2.}
\end{equation}

It follows from Floquet theory (see the next section) that 
$c\equiv r(H)>r(L)$. 
Hence, averaging over the eigenvalues of $L$ gives 
the following estimate (lower bounds):
\[
c^{2j}>{1\over q}\Tr L^{2j},\qquad j=1,2,\dots\label{simple-est}
\]

Below we obtain a different type of bounds for $c$:
\begin{theo}\label{Twidth}
{\rm (Bounds for the width of the spectrum)}  
Let $H$ be normalised and $q>1$. Then
\begin{eqnarray} 
&{}&c> 2A,\label{c}\\
&{}&c^2\left( {1\over2}+\ln{c\over 2A}\right)
>{1\over q}\Tr L^2.\label{c2}
\end{eqnarray}
\end{theo}

\no {\bf Remark:} The equality sign in (\ref{c}) happens if 
and only if $q=1$.

This theorem will be proved in Section 3.
The first part of it, (\ref{c}), will also 
be proved by a much simpler method which uses only general properties 
of polynomials (Lemma \ref{Lu}). Note that the quantity $A$ can be interpreted 
as the logarithmic capacity of the spectrum (see, e.g., \ci{VA}). 
One might ask if the inequality (\ref{c})
between half the diameter of a set and its logarithmic capacity holds 
in a more general context. 
For the Schr\"odinger case ($a_j\equiv 1$),
the inequality (\ref{c}) reduces to $c>2$. It was recently shown in \ci{KS}
that $c>2$ for an arbitrary, not necessarily periodic, real nonconstant 
sequence $b_j$, $a_j\equiv 1$.

As will be clear from the next section,
one can also derive further inequalities similar to (\ref{c2}). 
They are analogous to (\ref{simple-est}) but better. 
We obtain from them lower bounds for $c$. 
Consider, for example, Harper's matrix (see 
\ci{Harperrev} for a review):
$a_j\equiv 1$, $b_j=2\cos(2\pi\alpha j+\theta)$,
$\alpha=p/q$, $p$ and $q\in\Z$, $\theta\in\R$. 
Obviously, $c\le2+\max_j |b_j|\le 4$.
Simple calculation shows that $\Tr L^2=4q$.
Substituting this value into inequality (\ref{c2})
and solving the latter numerically gives
$c>2.41$.
Note that this bound does not depend on $p$ and $q$.
Therefore
it follows from continuity that it also holds for
irrational $\alpha$.

Let 
$$
h_+=\max_{\l\in \cup \g_n}{1\over q}{\rm arccosh}\;|A^{-q}\det(\l-L)/2|.
$$
We obtain the following estimates 
for the gaps of $H$ in terms of the matrix elements and $c$:

\begin{theo}\label{Tgaps}
{\rm (Bounds for the gaps)} 
Let $H$ be normalised and $q>1$. Furthermore,
let $g_n=(\arccos(\l^-_n/c),\arccos(\l^+_n/c))$, 
where $(\l^-_{n}, \l^+_n)$ is the $n$'th
gap in the spectrum of $H$. 
Then at least one gap is open and
\[
\sum_{n=1}^{q-1} |g_n|>\pi\frac{\ln(c/2A)}{h_+}
>\pi\frac{\ln(c/2A)}{\ln (2c/A)},\label{g1}
\]
\[
\sum_{n=1}^{q-1}\int_{g_n}\cos^2 x dx>
{\pi\over h_+}\left(\frac{\ln(c/2A)}{2}+{1\over4}-
{1\over 2qc^2}\Tr L^2\right),\label{g1p}
\]
\[
\frac{\ln(c/2A)}{\max\{1,qh_+/\pi\}}<
\sum_{n=1}^{q-1} |g_n|^2 < 8\ln(c/2A),\label{g2}
\]
where $h_+$ satisfies  
$0<h_+<\ln({c\over A}+|{1\over Aq}\sum_{j=1}^q b_j|)<\ln{2c\over A}$.
\end{theo}

\no {\bf Remarks:}
\begin{itemize}
\item Existence of an open gap implies that $h_+\neq 0$.
\item Just like for (\ref{c2}) it is possible to derive further 
inequalities similar to (\ref{g1p}) improving the estimates.
\end{itemize}


The plan of the paper is as follows.
In the next section we recall some facts from Floquet theory
and prove a lemma which guarantees existence of at least one open gap if 
$q>1$ and validity of (\ref{c}). Note that existence of gaps for the 
Schr\"odinger case also follows from (\ref{eL}) \ci{Last}.
In Section 3
we construct the quasimomentum, obtain the trace formulas and prove 
Theorem~\ref{Twidth}. In Section 4 we show how our 
construction fits into the theory of the general quasimomentum developed in
\ci{MO,L,KK1} and recall some properties of the general 
quasimomentum.
In Section 5 we use these properties together with estimates 
for the Dirichlet integrals (Lemma~\ref{Ldir}) to prove Theorem~\ref{Tgaps}
and to establish further bounds on the maxima of the Lyapunov 
exponent (Lemma~\ref{LQmax}). 
Note that the proof of the bounds (\ref{g1}) and (\ref{g1p}) is simpler
than that of (\ref{g2}). To establish the former, we only need 
Lemmas \ref{Lk}, \ref{Ltr}, and \ref{Lh+}.


\section{Preliminaries}
\setcounter{equation}{0}

Recall some facts from Floquet theory (see, e.g., \ci{Toda}).

For a fixed $\l$ we introduce two fundamental solutions 
$\phi(\l)$, $\theta(\l)$ of the equation ($\psi=\phi$ or $\psi=\theta$)
$(H\psi(\l))_n=\l\psi_n(\l)$, $n=-\infty,\dots,\infty$, 
by the initial conditions
$\phi_0=0$, $\phi_1=1$; $\theta_0=1$, $\theta_1=0$.
Both of them satisfy the recurrence relation:
\[
a_n\psi_{n+1}(\l)=(\l-b_n)
\psi_n(\l)-a_{n-1}\psi_{n-1}(\l), \qquad n=1,2,\dots.\label{rec}
\]
The polynomial of degree $q$
\[
D(\l)=\phi_{q+1}(\l)+\theta_q(\l)\label{discriminant}
\]
is called the discriminant of $H$.

For example, when $q=2$ we have $D(\l)
=({\l}^2-(b_1+b_2){\l}+b_1b_2-a_1^2-a_2^2)/(a_1a_2)$.

It is known that:

\no 1) Any solution of the equation $(H-\l)\psi=0$ 
(except possibly for those corresponding to $D(\l)=\pm 2$)
is a linear combination of two solutions $\psi^+$ and $\psi^-$ with 
the property
\[
\p^\pm_{n+q}(\l)=e^{\pm iq\wt k(\l )}\p^\pm_n(\l),  \quad n\in \Z,
\quad {\rm where } \quad  \wt k(\l )={1\over q}
\arccos {D(\l )\over 2}.\label{p1}
\]
Henceforth, we shall fix the branch of $\arccos(x)$ by the condition
$\arccos(0)=-\pi/2$.
The function $\wt k(\l)$ plays a crucial role in the present paper. 

\no 2) The spectral intervals $\s_m$ are 
the image of $[-2,2]$ under the inverse of the
transform $D(\l)$ (see Figure 1). (Note that by (\ref{p1}) solutions
$\psi_n$ are bounded in $n$ if $D(\l)\in [-2,2]$ and 
otherwise exponentially increase.) Hence,  
all the critical points of $D(\l)$ are
maxima and minima; moreover, at all maxima points $\l_{\rm max}$
$D(\l_{\rm max})\ge2$ and at all minima points $\l_{\rm min}$
$D(\l_{\rm min})\le-2$. All $q-1$ critical points are mutually different.

\no 3) The solutions $\phi(\l)$ and $\theta(\l)$ satisfy the identity
\[
\phi_{q+1}(\l)\theta_q(\l)-\phi_{q}(\l)\theta_{q+1}(\l)=1.\label{Wron}
\]

\begin{figure}
\centerline{\psfig{file=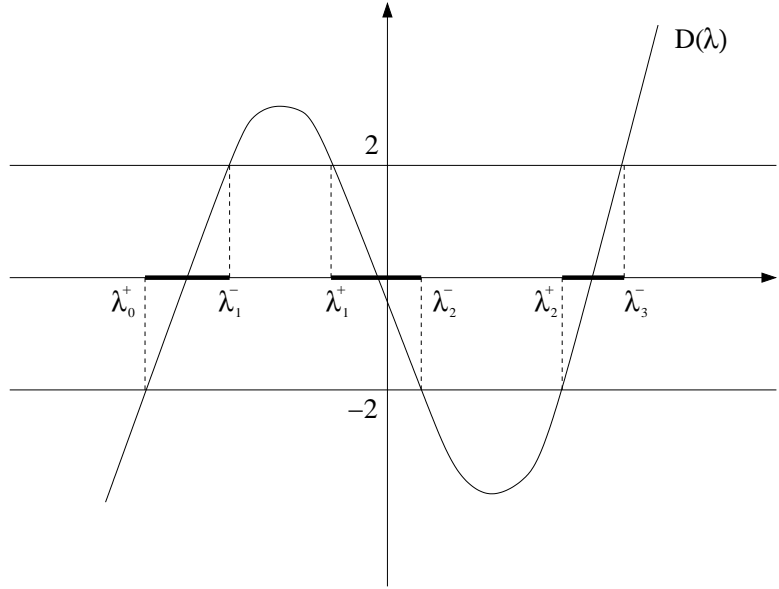,width=4.5in,angle=-90}}
\vspace{0.1cm}
\caption{
The discriminant is sketched for $q=3$.
The bands of the spectrum are shown by thick lines.}
\label{fig1}
\end{figure}

We shall need the representation of $D(\l)$ given by

\begin{lemma}\label{LD}
$D(\l)=A^{-q}\det(\l-L)$ for all $\l\in\C$.
\end{lemma}

\no {\it Proof.} Denote $D_{jk}$, $j,k>0$, the determinant of the 
matrix $\l-L$ with
the first $j-1$ and the last $q-k$ rows and columns removed, and
$L_{1q}$, $L_{q1}$ set to zero. 
Expanding $\det(\l-L)$ in the elements of the first row, we get 
\[
\det(\l-L)=D_{1q}-a_q^2 D_{2\,q-1}.\label{det}
\]

Similarly, the expansion of $D_{1j}$ by the last row gives: 
$D_{1j}=(\l-b_j)D_{1\,j-1}-a_{j-1}^2D_{1\,j-2}$, $j=1,2,\dots,q$,
where $D_{1\,-1}=0$, $D_{10}=1$.
Comparing this recurrence with (\ref{rec}), we get
$D_{1q}=a_1a_2\cdots a_q\phi_{q+1}$ by induction.
Similarly, we show that 
$D_{2\,q-1}={a_1a_2\cdots a_{q-1}\over -a_q}\theta_q$. Substituting
these expressions into (\ref{det}) proves the lemma. $\BBox$ 

Now we shall prove that $q\neq 1$ implies the inequality
(\ref{c}) and existence of gaps. 
The proof is based on the above given facts from Floquet theory 
and general properties of polynomials.

\begin{lemma}\label{Lu}
Let $q>1$. Then $c=r(H)>2A$ and at least one gap is open.
\end{lemma}

\no {\it Proof.}
Suppose that $c\le 2A$. 
Then by adding a constant to $H$, we first ensure 
that the spectrum of $H$ (${\rm Spec}(H)$) lies within the interval 
$I=[-2A,2A]={\rm Spec}(H_0)$. Here $H_0$ is the Jacobi matrix
with $a_j\equiv A$ and $b_j\equiv 0$. 
Denote by $D_H(\l)$ the discriminant of $H$.
It is a polynomial of degree $q$ with the coefficient $A^{-q}$ of $\l^q$.
Let us set 
all $a_j=A$ and $b_j=0$ in it, and denote the resulting 
polynomial $D_{H_0}(\l)$. 
It is the discriminant for $H_0$ viewed as a matrix of period $q$.
Since 
${\rm Spec}(H_0)$ has no gaps, we conclude that 
$\max_{\l\in I}|D_{H_0}(\l)|=2$. 

Let $\l_0=-2A$, $\l_q=2A$, and $\l_j$, $j=1,\dots,q-1$ be the 
critical points of 
$D_{H}(\l)$. By our construction, they all lie in $I$. By Floquet theory, 
$|D_H(\l_j)|\ge2$ and the signs in the sequence $D_H(\l_j)$, $j=0,\dots,q$ 
alternate. Thus
\[
D_H(\l_q)\ge D_{H_0}(\l_q),\quad 
D_H(\l_{q-1})\le D_{H_0}(\l_{q-1}),\quad 
D_H(\l_{q-2})\ge D_{H_0}(\l_{q-2}),\quad\dots\label{ineqD}
\] 
The polynomial $S(\l)=D_H(\l)-D_{H_0}(\l)$ is of degree less than $q$.
However, as follows from (\ref{ineqD}), it changes its sign
on $I$ at least $q$ times. Therefore $S(\l)\equiv 0$.
Thus the strict inequality $c<2A$ is impossible, and $c=2A$  
happens only when $D_H\equiv D_{H_0}$. We shall now demonstrate that the 
last identity implies $H=H_0$, i.e., the smallest period of $H$ is 1. 

The method is essentially borrowed from \ci{Perkolab}.
Assume $D_H\equiv D_{H_0}$ (we shall omit the subscript).
Hence $\l^+_j=\l^-_j$ and $\max_{\l\in I}|D(\l)|=|D(\l^\pm_j)|=2$.
Let $\nu_j$, $j=1,\dots,q-1$ be the zeros of the $(q-1)$-degree polynomial 
$\phi_q(\l)$. Since $\phi_n(\l)$, $n=1,2,\dots$ are the 
orthogonal polynomials corresponding to $H$,
all the $\nu_j$ are simple and belong to the open set $(-2A,2A)$.
Using (\ref{Wron}), we get
$|D(\nu_j)|=|\theta_q(\nu_j)^{-1}+\theta_q(\nu_j)|\ge 2$.
Therefore $\nu_j=\l^\pm_j$, $|D(\nu_j)|=2$, and $\theta_q(\nu_j)=\pm 1$,
where the sign corresponds to that of $D(\nu_j)$, $j=1,\dots,q-1$.
Since $\theta_q(\l)$ is a polynomial of degree $q-2$, it can be uniquely 
reconstructed from its values in $q-1$ points $\nu_j$ by the Lagrange 
interpolation formula.
From the definition of the discriminant, we obtain 
$\phi_{q+1}(\l)=D(\l)-\theta_q(\l)$.
On the other hand, we also know the monic polynomial (i.e., with the 
coefficient 1 of the highest degree) 
$\hat\phi_q(\l)=\prod_{j=1}^{q-1}(\l-\nu_j)$; note that it is $D'(\l)$ 
up to a factor.
We shall now see that these two polynomials determine $H$.
Note first that the coefficient of the highest ($n-1$) degree of $\phi_n(\l)$
is $(a_1a_2\cdots a_{n-1})^{-1}$. Therefore we have the following recursion
for the monic polynomials 
$\hat\phi_n(\l)=\l^{n-1}+\alpha_n^{n-2}\l^{n-2}+\cdots+\alpha_n^0$:
\[
\hat\phi_{n+1}(\l)+(b_n-\l)\hat\phi_{n}(\l)+a^2_{n-1}\hat\phi_{n-1}(\l)=0.
\label{rec1}
\]
All coefficients of this polynomial are zero. This gives
\[
b_n=\alpha_n^{n-2}-\alpha_{n+1}^{n-1}, \qquad
a_{n-1}=\sqrt{\alpha_n^{n-3}-\alpha_{n+1}^{n-2}-b_n\alpha_n^{n-2}}.
\]
Assuming $\hat\phi_{n+1}(\l)$ and $\hat\phi_{n}(\l)$ are given,
we then obtain the polynomial $\hat\phi_{n-1}(\l)$ from (\ref{rec1}).
Performing this procedure successively for $n=q,q-1,\dots,1$, we 
reconstruct the coefficients $a_1,\dots,a_{q-1}$ and $b_1,\dots,b_q$.
Finally, $a_q$ is obtained from the highest-degree coefficient $\alpha$ of 
$\phi_{q+1}(\l)$:   $a_q=(\alpha a_1\cdots a_{q-1})^{-1}$.
From the uniqueness of our reconstruction, it follows that $H=H_0$. 

Thus for $q\ge 2$ we have $c>2A$. Suppose now that the spectrum of $H$ has 
no gaps, i.e., $|D_H(\l)|\le 2$ on $[-c,c]$ (assuming $H$ is normalised).
Then considering the inequalities between $D_H(\l)$ and $D_{H_0}(\l)$
now at the critical points of $D_{H_0}(\l)$ on $[-c,c]$, we obtain as before
$S(\l)=D_H(\l)-D_{H_0}(\l)\equiv 0$. This contradiction shows that there 
is at least one open gap in the spectrum of $H$. 
$\BBox$

\section{Trace formulas and Dirichlet integrals}
\setcounter{equation}{0}

Our first aim is to investigate the asymptotics of the function 
$\wt k(\l)$ for large 
imaginary $\l$. We shall see below that it is convenient to introduce 
instead of $\l$ a new variable $z$ such that the mapping $\wt k(\l(z))$
becomes asymptotically the identity. 
After determining the asymptotics of $\wt k(\l(z))$, we construct 
the regions into 
which the upper half-plane is mapped by the functions $z(\l)$ and 
$\wt k(\l)$.
This information is then used to obtain the trace formulas and 
Dirichlet integrals for the quasimomentum.


First, we set $c=\l^-_q=|\l^+_0|$ (recall the normalisation 
of section 1). 
Then, in the variable $\zeta=\l/c$, the spectrum lies within the 
interval $[-1,1]$ and its boundaries coincide with $-1$ and $1$.

Using the expansion
\[
\arccos x=i\left(\ln 2x-\sum_{j=1}^\infty\pmatrix{2j\cr j}{1\over 2j
  (2x)^{2j}}\right), \qquad |x|>1,\label{arccos} 
\]
we get $q\wt k(\l)=\arccos(D(\l)/2)=i\ln D(\l)+ O(1/\l^{2q})$
as $|\l|\to\infty$.
Here, by Lemma~\ref{LD},
$$
\eqalign{
\ln D(\l)=\ln\det(\l-L)-\ln A^q=\Tr\ln(\l-L)-\ln A^q=\\
q\ln{\l\over A}+\Tr\ln(1-L/\l)= 
q\ln{\l\over A}-\sum_{j=1}^\infty\frac{\Tr L^j}{j\l^j}.}
$$

Further, using (\ref{arccos}), we have for $|\zeta|>1$
$$
i\ln{\l\over A}=i\ln{c\over 2A}+i\ln 2\zeta=
i\ln{c\over 2A}+\arccos \zeta +i\sum_{j=1}^\infty\pmatrix{2j\cr j}{1\over 2j
  (2\zeta)^{2j}}.
$$

\no {\bf Definition.} The function
\[
k(z)=\wt k(-c\cos z)+\pi,
\]
where
$\wt k(\l)$ is defined by (\ref{p1}), (\ref{discriminant}),
is called quasimomentum. 

From the preceding discussion, we have:
\[
k(z)=z+iQ_0+\sum_{j=1}^\infty{iQ_j\/\cos^j z},\qquad z\to i\iy,\label{asymp}
\]
where
\[
Q_0=\ln{c\over 2A};\qquad
Q_j=\cases{
(1/ jc^j){1\over q}\Tr L^j & if $1\le j<2q$ and odd;\cr
(1/ j2^j)\pmatrix{j\cr j/2}-
(1/ jc^j){1\over q}\Tr L^j  & if $2\le j<2q$ and even.}\label{Q}
\]
In particular, 
\[
Q_1={1\over qc}\Tr L,\quad
Q_2={1\over 4}-{1\over 2qc^2}\Tr L^2.
\]


Investigating the mapping of the boundaries (as we explain below), 
we obtain the correspondence 
of the domains in $\l$-, $z$-, and $k$-planes under the conformal 
transformations $z(\l)=\pi+\arccos\l/c$, $k=k(z)$ as shown 
in Figures 2--4 (for $q=3$).
The notation for the real and imaginary parts of $z$ and $k$ 
is fixed as follows: $z=x+iy$, $k=u+iv$. 
The points $\l^+_n$, $z^+_n$, $k^+_n$ correspond to each other under 
the mappings;
so do the points $\l^{-}_n$, $z^{-}_n$, $k^{-}_n$.
That the upper $\l$-half-plane is mapped onto the half-strip by the
function $z(\l)$ is easy to verify. We prove first the mapping of the
boundary $\partial S^\l_R$ onto $\partial S_R$ (for finite regions
$S^\l_R$ and $S_R$) and then let $R\to\infty$.
To obtain the mapping of the boundary of $\l$ domain onto that of $k$
domain, it is convenient to represent $\wt k(\l) +\pi$ in the form
$$
{1\over q}\int^{D(\l)}_{D(\l_0^+)}\frac{dD(\l)}{\sqrt{4-D(\l)^2}}
$$
and integrate (keeping in mind Figure 1) along 
the real part of $\partial S^\l_R$ with
vanishing semicircles above the ends of the gaps $\l^\pm_j$.

\begin{figure}
\centerline{\psfig{file=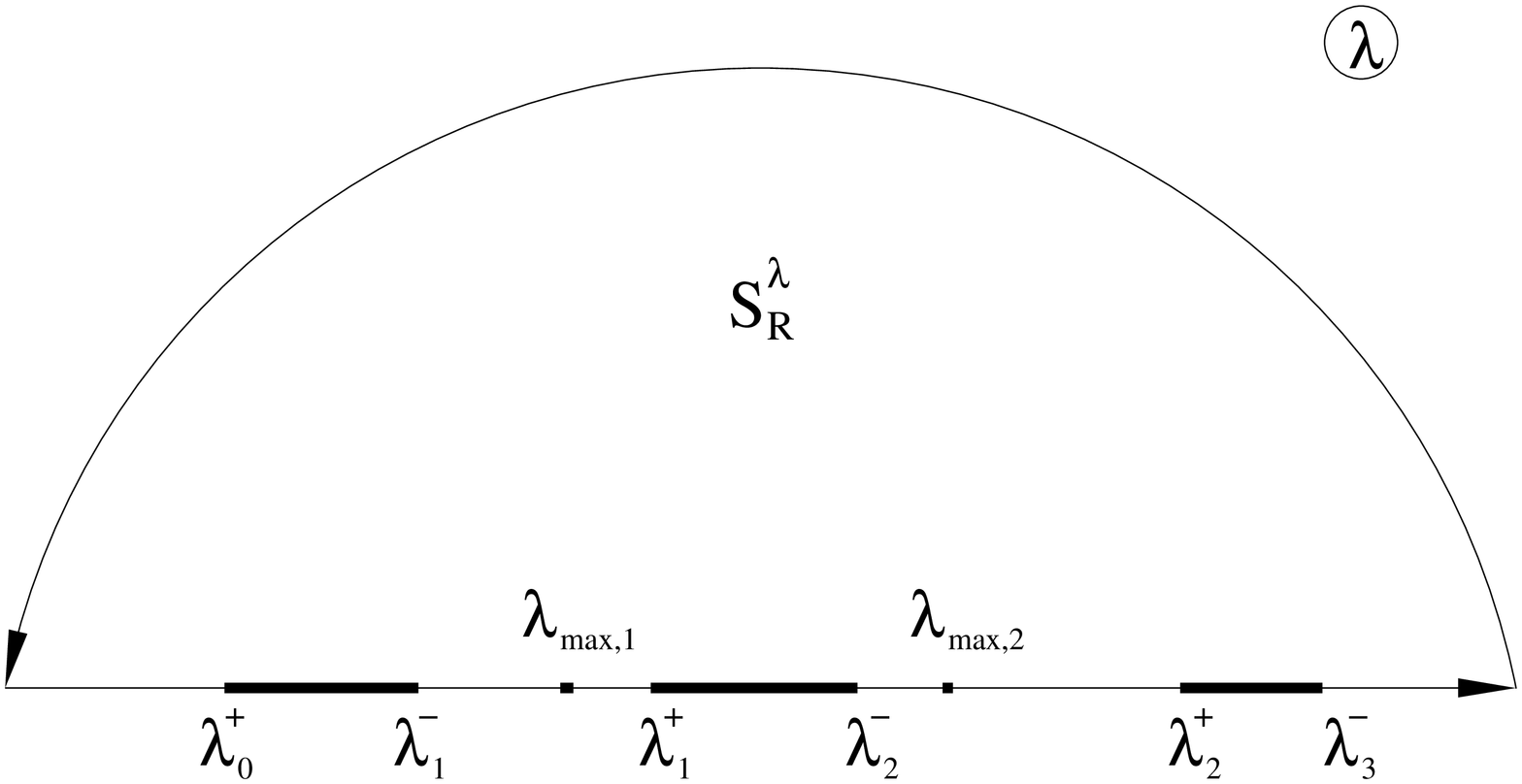,width=4.5in,angle=-90}}
\vspace{0.1cm}
\caption{
$\l$-plane for $q=3$.}
\label{fig2}
\end{figure}
\begin{figure}
\centerline{\psfig{file=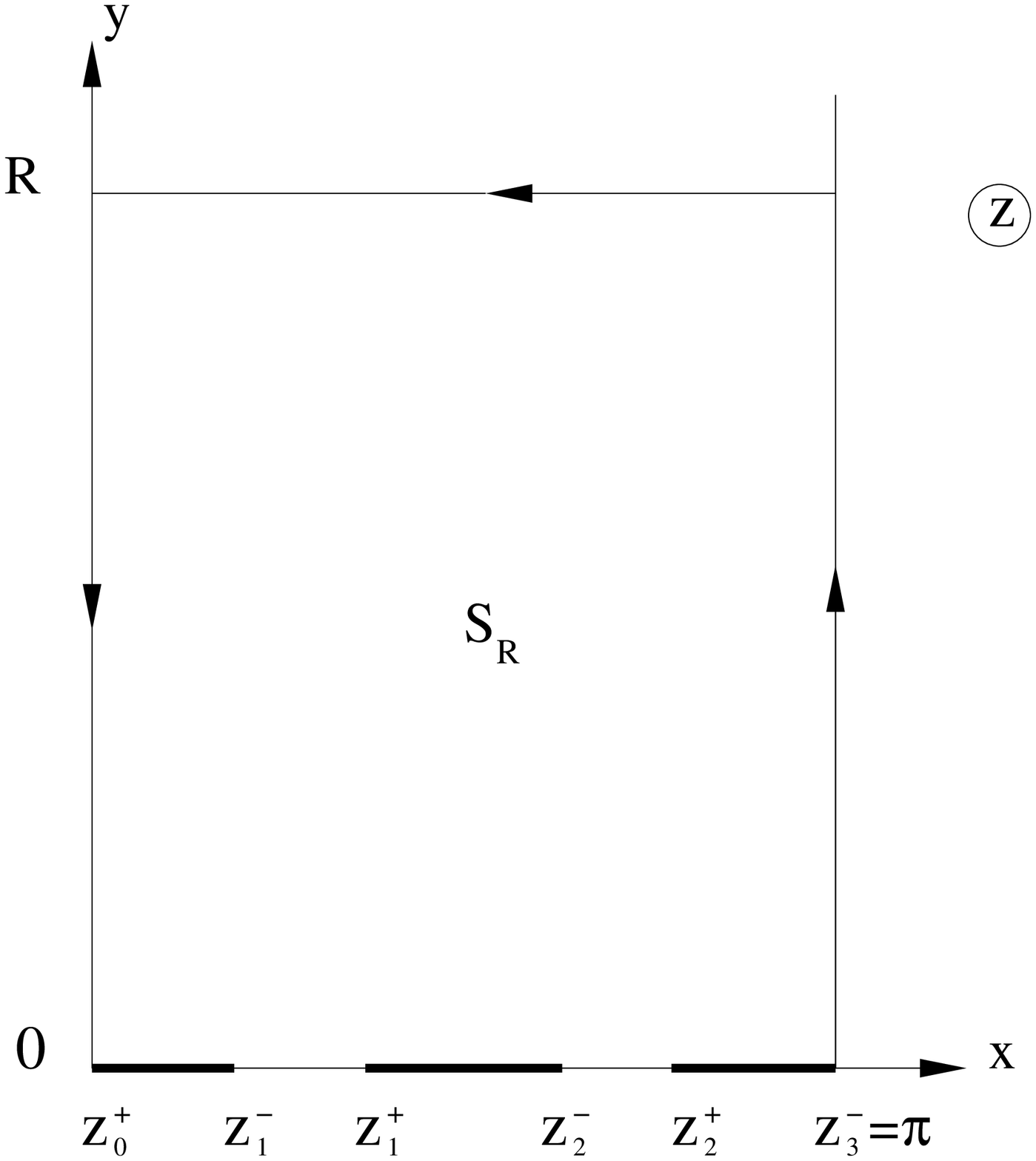,width=4.5in,angle=-90}}
\vspace{0.1cm}
\caption{
$z$-plane.}
\label{fig3}
\end{figure}
\begin{figure}
\centerline{\psfig{file=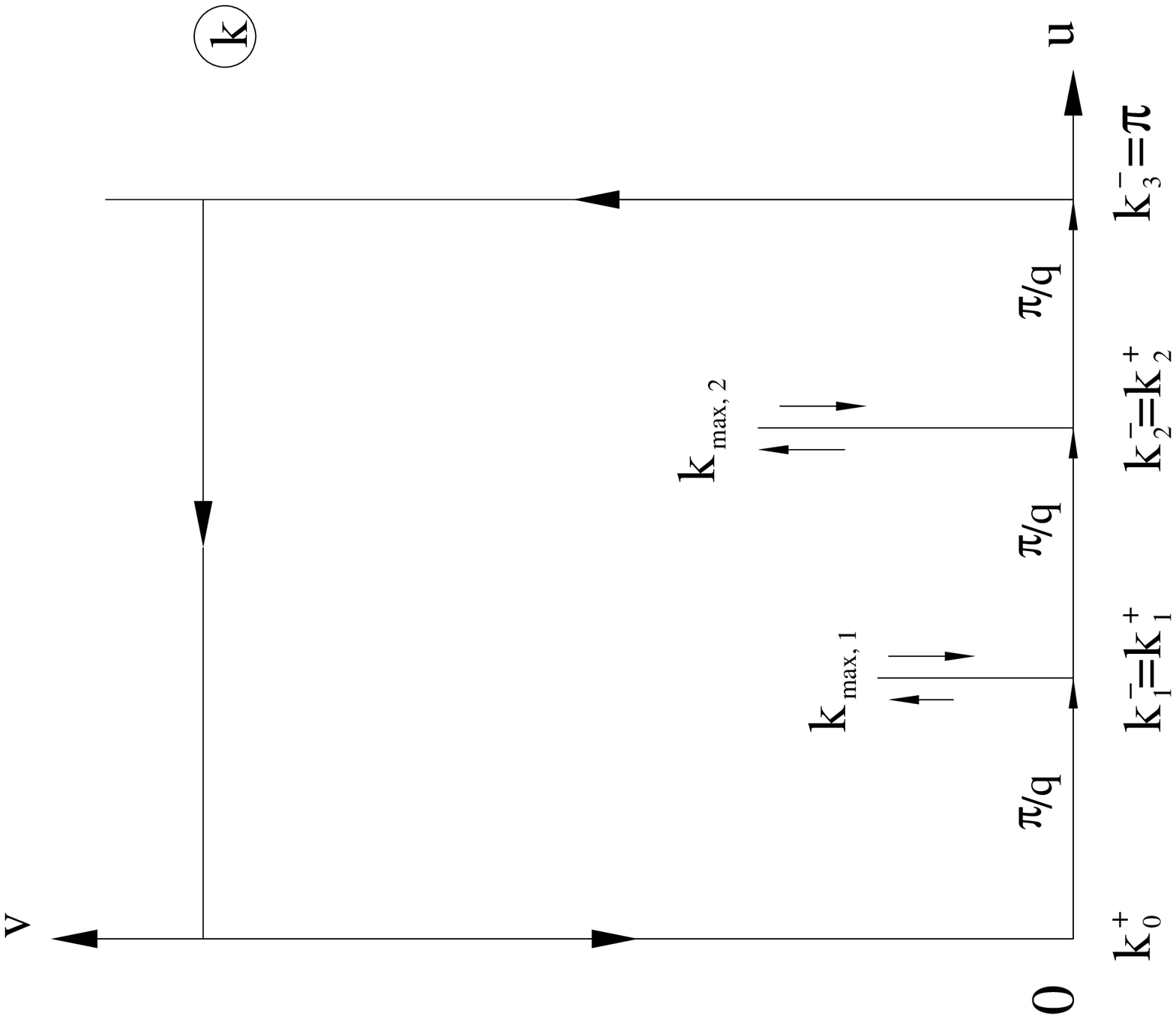,width=4.5in,angle=-90}}
\vspace{0.1cm}
\caption{
$k$-plane}
\label{fig4}
\end{figure}

We now write these results more precisely. Let
\begin{eqnarray}
\Lambda_+&=&\ol\C_+\sm \cup \g_n,\qquad \C_+=\{\l:\Im \l>0\},\nonumber
\\
\cZ_+^0&=&\C_+^0\sm  \cup g_n ,\qquad \C_+^0=\{z: 0\le\Re z\le\pi, 
\Im z\ge0\},
\nonumber
\\
g_n&=&(z^{-}_n,z^+_n),\qquad 
n=1,\dots,q-1,\nonumber
\end{eqnarray}
where $z^\pm_n=\pi+\arccos \l^\pm_n/c$. 
The function $z=\pi+\arccos \l/c$
is a conformal mapping of $\Lambda_+$ onto $\cZ_+^0$. The gap $\g_n$ 
is mapped onto $g_n$, $n=1,\dots,q-1$. 
Let
$$
\cK_+^0=\C_+^0\sm \cup c_n,\qquad 
 c_n=({n\pi \/q},{n\pi \/q}+ih_n),\qquad n=1,\dots,q-1,
$$
where $c_n$ is a vertical slit.
The function $k=k(z)$ conformally maps $\cZ_+^0$ onto $\cK_+^0$. 
The spectral interval 
$[z^+_{n-1},z^{-}_n]$ is mapped onto the interval of length $\pi/q$ of the 
$u$-axis $[(n-1)\pi/q,n\pi/q]$, $n=1,\dots,q$; 
the gap $g_n$, onto the slit $c_n$, $n=1,\dots,q-1$.
The height of the $n$-th slit is $h_n=(1/q){\rm arccosh}\;
|D(\l_{{\rm max},\, n})/2|$, where
$\l_{{\rm max},\,n}$ is the critical point of $D(\l)$ in the gap $\g_n$.

Thus we proved

\begin{lemma}\label{Lk}{\rm (quasimomentum)}
The function $k(z)=\wt k(-c\cos z)+\pi$ is a conformal mapping of
the domain $\cZ_+^0$ onto the quasimomentum domain $\cK_+^0$. 
It possesses the asymptotic expansion (\ref{asymp}).
\end{lemma}

We shall now obtain the identities (trace formulas) which connect the
integrals of $u\equiv\Re k$ and $v\equiv\Im k$ with the coefficients 
$Q_i$ in the asymptotic 
expansion of $k(z)$ and, via these coefficients, with the matrix 
elements of $H$. 
The trace formulas will be the basis of all the spectral estimates in 
the present work. 
To establish (\ref{c}) (in a way different from that of Lemma~\ref{Lu}), 
only the first of these formulas is sufficient.
The other ones can also be used to extract some local information about 
the gaps, but we shall not do this here.

\begin{lemma}\label{Ltr}{\rm (trace formulas)} 
The following identities hold:
\begin{eqnarray}
{1\over\pi}\int_0^\pi v(x)dx&=& Q_0=\ln{c\/2A},\label{v0}\\
{1\over\pi}\int_0^\pi v(x)\cos^n xdx&=&
\cases{\sum_{i=0}^{(n-1)/2}Q_{2i+1}\pmatrix{n-1-2i\cr (n-1)/2-i}
/2^{n-1-2i}& if n is odd,\cr
\sum_{i=0}^{n/2}Q_{2i}\pmatrix{n-2i\cr n/2-i}
/2^{n-2i}& if n is even.}\label{vall}
\end{eqnarray}
In particular,
\begin{eqnarray}
{1\over\pi}\int_0^\pi v(x)\cos x dx&=&Q_1,\label{v1}\\
{1\over\pi}\int_0^\pi v(x)\cos^2 x dx&=&Q_0/2+Q_2.\label{v2}
\end{eqnarray}
\end{lemma} 

\no {\it Proof.}
Let us calculate the integral of $k(z)$ along the boundary of $\cZ_+^0$ as 
shown in Figure 3. We integrate first along the contour $\partial S_R$:
$$
0=\int_{\partial S_R}k(z)dz=\int_0^\pi k(x)dx+
\int_0^R(\pi+iv(\pi,y))idy+
\int_\pi^0 k(x+iR)dx+
\int_R^0 iv(0,y)idy.
$$
Now take the limit $R\to\infty$ using (\ref{asymp}) for $k(x+iR)$.
We obtain
$$
0=\int_{\partial\cZ_+^0}k(z)dz=\int_0^\pi k(x)dx-i\pi Q_0-\pi^2/2
-\int_0^\infty(v(\pi,y)-v(0,y))dy.
$$
Separating the real and imaginary parts, we obtain (\ref{v0}) and 
the identity
\[
\int_0^\pi u(x) dx={\pi^2\over2}+
\int_0^\infty(v(\pi,y)-v(0,y)) dy.\label{u}
\]

Similarly, considering the integrals of
$(k(z)-z-iQ_0-i\sum_{j=1}^nQ_j/\cos^j z)\cos^n z$, $n=1,2,\dots$ along
$\partial\cZ_+^0$ we get
\[
\int_0^\pi v(x)\cos^n xdx=
\sum_{j=0}^n Q_j\int_0^\pi\cos^{n-j}xdx,
\]
whence the lemma follows. $\BBox$


{\it Proof of Theorem~\ref{Twidth}.}  

i) Since $v(x)\ge 0$, it follows from (\ref{v0})
that $c\ge2A$. 

ii) Similarly, the right-hand side of (\ref{v2}) 
is nonnegative. Substituting there 
(\ref{Q}), we obtain (\ref{c2}) with the greater-or-equal sign. 

To make the inequalities strict, we use Lemma~\ref{Lu}: if $q>1$ there exists 
at least one open gap, hence $v(x)$ is not identically zero. Thus the 
integrals in Lemma~\ref{Ltr} are strictly positive and the theorem is proven.
$\BBox$

Similarly, using (\ref{vall}), we can obtain further inequalities for $c$. 

For a function $f(z), z\in \C_+$, we formally define the Dirichlet integral
$$
I(f)={1\/\pi}\iint_{\cZ_+^0}|f'(z)|^2dxdy.
$$
\no We shall now present two such integrals the first of which 
we use below for our estimates.

\begin{lemma}\label{Ldir}{\rm (Dirichlet integrals)} 
The following identities hold:
\begin{eqnarray}
I(k(z)-z)&=&Q_0=\ln{c\/2A};\label{Dir1}\\
I\left(\left[k(z)-z-iQ_0\right]\cos z\right)&=&
{Q_0\over2}+Q_2-2Q_0Q_2-{Q_1^2\over2}.\label{Dir2}
\end{eqnarray}
\end{lemma}

\no {\it Proof.} 
Using the Cauchy conditions, we rewrite 
the integral $I(f(z))$ in the form
\[
I(f(z))={1\over\pi}\iint_{\cZ_+^0}|\nabla\Im f(z)|^2dxdy.
\]
Recall the Green formula for the harmonic function 
$\omega(x,y)=\Im f(z)$ in the closed area $S_R$:
\[
\iint_{S_R}|\nabla\omega(x,y)|^2dxdy=\int_{\partial S_R}\omega
\frac{\pa \o }{\pa n}dl,
\label{Gauss}
\]
where $\pa /{\pa n}$ denotes the normal derivative
(the normal points outside $S_R$) to the 
curve $\partial S_R$. 

Let us calculate (\ref{Gauss}) for $f=k(z)-z$, that is for 
$\omega(x,y)=v(x,y)-y$.
We use the Cauchy conditions to see that the integrals along both vertical
boundaries of $S_R$ vanish.  
To calculate the integral along the interval $[0,\pi]$ of the real
line, we note
(cf. Figure 4) that $u={\rm const}$ on each gap, and  $v=0$  on each
band, hence the term $v(x)\partial u(x)/\partial x$ vanishes on
$[0,\pi]$ and the contribution of this part of the path is
$\int_0^\pi v(x)dx$.
Taking then the limit as $R\to\infty$ (using the asymptotics of
$k(z)$ on the part of $\partial S_R$ where $y=R$ to see that the
contribution of that part to the integral tends to zero), we obtain
\[
\iint_{\cZ_+^0}|k'(z)-1|^2dxdy=\int_0^\pi v(x)dx.  
\]
In view of (\ref{v0}), we have (\ref{Dir1}).

We turn to the calculation of $I(f)$ for
$f=\left[k(z)-z-iQ_0\right]\cos z$.
Proceeding as above and using (\ref{v2}), we get
\begin{eqnarray}
I\left(\left[k(z)-z-iQ_0\right]\cos z\right)=
{3\over4}Q_0+Q_2+{1\over2\pi}I',\label{I2mid}\\
I'=\int_0^\pi(vu-vx+Q_0 u)\sin 2x dx.\label{Ip}
\end{eqnarray}
In order to evaluate $I'$, we consider
\begin{eqnarray}
0=\lim_{\epsilon\to 0}\int_{\partial_\epsilon\cZ_+^0}
\left({k(z)^2\over2}-k(z)z+iQ_0k(z)+{z^2\over2}+{3\over2}Q_0^2-iQ_0z+
\right.
\nonumber\\
\left.
{2Q_0Q_1\over\cos z}+{2Q_0Q_2+Q_1^2/2\over\cos^2 z}\right)\sin2z dz,
\label{addi}
\end{eqnarray}
where $\partial_\epsilon\cZ_+^0$ differs from $\partial\cZ_+^0$ only in the 
$\epsilon-$neighbourhood of the point $z=\pi/2$ where it is a semicircle above
this point of radius $\epsilon$ (the integrand has a simple pole at $z=\pi/2$).
The integrand in (\ref{addi}) is chosen to satisfy two conditions:
1) integration of the imaginary part of the first three terms in it 
over $[0,\pi]$ gives
$I'$ (\ref{Ip}); 2) the integrand behaves like $O(1/\cos z)$ as $z\to i\infty$.
We evaluate (\ref{addi}) as the contour integrals in the proof of 
Lemma~\ref{Ltr}. The imaginary part of (\ref{addi}) then yields
$$
I'=-\pi\left({Q_0\over2}+4Q_0Q_2+Q_1^2\right).
$$
Substituting this into (\ref{I2mid}), we conclude the proof of the lemma.
$\BBox$


Further formulas similar to (\ref{Dir1}, \ref{Dir2}) can be obtained 
with more effort.

\section{Properties of the quasimomentum}
\setcounter{equation}{0}

By Lemma~\ref{Lk}, the quasimomentum $k(z)$ conformally maps the strip
$\cZ_+^0$ onto the strip $\cK_+^0$. We can expand it to a 
conformal mapping of the upper half $z$-plane onto the upper half $k$-plane 
(with slits) in the following way.
First, let $k(-x+iy)=k(x+iy)$ for $x+iy\in\cZ_+^0$ (reflection 
with respect to the $y$ axis). Now continue $k(z)$ periodically:
$k(z+2\pi n)=k(z)$, $n=\pm 1,\pm 2,\dots$.
Thus extended $k(z)$ is a particular (periodic) case of the general 
quasimomentum: a well-studied mapping \ci{L,KK1}. 
We can utilise, therefore, the known properties 
of the general quasimomentum. The rest of this section will be devoted to a
review of some of these properties.  

Define the so-called comb domain
$$
\cK_+=\C_+\sm \G, \ \ \G=\cup \G_n,\ \ \  \ \G_n=(u_n,u_n+ih_n),\ \  
h_n \geq 0,\ \        n\in\Z,
 $$
where  $u_n$ is a strongly increasing sequence of real numbers
such that $u_n\to \pm\iy$ as $n\to\pm\iy$, and
$\{h_n\}_{-\iy}^{\iy}\in \ell^{\iy}$.

A conformal  mapping $k(z)=u+iv$ from the upper
half  plane $\C_+$ onto some comb $\cK_+$ is called a general
quasimomentum 
if  $k(0)=0$,  and $k(iy)=iy(1+{\rm o}(1))$ as $y\to \iy$. 
A general quasimomentum  $k(z)$
is a continuous function in $z\in \ol \C_+$. 
The inverse function $z(k)$ maps each slit $\G_n$ onto the interval (gap) 
$g_n=(z^-_n,z^+_n)$ of the real axis; and each interval $[u^+_{n-1},u^-_n]$
onto (the band) $[z^+_{n-1}, z^-_n]$.
Obviously, $v(z)>0$ for $z\in\C_+$.
The function $v(x)$ is continuous on the real axis, equals zero on the 
bands, is positive and reaches the maximum $h_n$ in the $n$'s gap.
Hence $v\in L^{\iy}(\R)$.
It is known that also
$u'_x(z)=v'_y(z)>0$ for $z\in\C_+$. On the real line 
the function $u'(x)$ 
is positive on the bands and equals zero on the gaps. (Thus it plays
a similar role for the bands as $v(x)$ for the gaps. However, it 
is more difficult to obtain good estimates for $u'(x)$.)
By the Herglotz theorem for positive harmonic functions (e.g., 
\ci{Koosis}), 
we have
\[
v(z)=y\left(1+{1\over \pi}\int_{-\infty}^\infty 
{v(t)\/ |t-z|^2}dt\right),\ \ \ z\in \C_+.\label{Herglotz_v}
\]
For its harmonic conjugate:
\[
u(z)=x+{1\over \pi}\int_{-\infty}^\infty 
\left({t-x \over |t-z|^2}-{t\over t^2+1}\right)v(t)dt+\const,\ 
\ \ z\in \C_+.\label{Herglotz_u}
\]
(The term with ${t\over t^2+1}$ is here only to ensure convergence if
$v(t)$ does not decay at infinity.)

Further,
\[
k(z)=u+iv=z+C_0+
{1\over \pi}\int_{-\infty}^\infty v(t)\lt({1\/ t-z}-{t\/ 1+t^2}\rt)dt,
\ \ \ \ C_0=-{1\/ \pi}\int_{-\infty}^\infty 
{v(t)dt\/ t(1+t^2)},\label{Herglotz_k}
\]
where the value of the constant $C_0$ 
is obtained from the condition $k(0)=0$.

Since by the definition of the quasimomentum 
$v'_y(iy)=1+o(1)$ as  $y\to\iy$,
we have $u'_x(iy)=1+o(1)$ as $y\to\iy$.
Hence, the Herglotz representation for $u'_x$ has the form
\[
u'_x(z)=1+{1\over \pi}\int_{-\infty}^\infty 
{y u'(t)\/ |t-z|^2}dt,\ \ \ z\in \C_+.\label{Herglotz_up}
\]
Consequently, on each gap $g_n$,
\[
-v''_{xx}|_{y\downarrow 0}=(u'_x)'_y|_{y\downarrow 0}=
{1\/\pi}\int_{-\infty}^\infty {u'(t)dt\/(t-x)^2}>0, \ \ \ x\in g_n,
\]
i.e. the function $v(x)$ is {\it concave} on each gap.

For our estimates we shall need the following result
proved in \ci{KK1}:
For any gap $g_j$ one has
\[
v(x)=w_j(x)\lt(1+{1\/ \pi}\int_{\R \sm g_j}{v(t)dt\/|t-x|w_j(t)}\rt),\ \ 
w_j(z)=|(z-z^+_j)(z-z^-_j)|^{1/2},\ \ \  \ x\in g_j.
\]

\no Thus, the graph of the function $v(x)$ 
lies {\it above the semicircle} over the gap (see Figure 5).

\begin{figure}
\centerline{\psfig{file=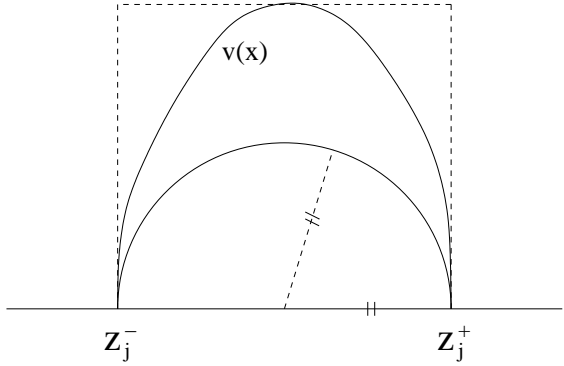,width=3.2in,angle=-90}}
\vspace{0.1cm}
\caption{
The graph of $v(x)$ is sketched on a gap $x\in(z_j^{-},z_j^+)$.}
\label{fig5}
\end{figure}

Note that in our case the quasimomentum is periodic (with period $2\pi$).
This allows us to reduce the general formulas. 
For example, (\ref{Herglotz_k}) takes on the form

\begin{lemma}\label{Lkp}
Let $k(z)=k(z+2\pi)$, $z\in \C_+$. Then
\[
k(z)=z+C+{1\over 2\pi}\int^\pi_{-\pi}v(\theta)
{\rm \; cotan}\left(
{\theta-z\over 2}\right)d\theta,\qquad 
C=-{1\over 2\pi}{\rm V.p.}\int^\pi_{-\pi}v(\theta){\rm \; cotan}
{\theta\over 2}d\theta.\label{Herglotz_kmod}
\]
\end{lemma}

\no{\it Proof.} We represent (\ref{Herglotz_k}) as a sum of integrals over 
the intervals $[-\pi+2\pi n, \pi+2\pi n]$, $n=-\infty,...,\infty$,
change the integration variable in each of them to transform it to 
$[-\pi,\pi]$, and then use the periodicity of $v(x)$ and the identity
$$
{\rm V.p.}\sum_{n=-\infty}^\infty{1\over \theta-z+2\pi n}={1\over 2}
{\rm\;cotan}\left({\theta-z\over 2}\right).
$$
The value of the constant is obtained from the condition $k(0)=0$.
$\BBox$

\section{Estimates}
\setcounter{equation}{0}
We return to the quasimomentum we constructed in Section 3.
Recall that for $q>1$ the function $v(x)$ is not identically zero.
Recall also that we assume $H$ to be normalised.

First, as a side result, note that
similarly to the case of the Hill equation \ci{KK1}, we have 
a double-sided estimate for $Q_0$ in terms of the gap lengths $|g_n|$ 
and maxima $h_n$ of $v(x)$ in the gaps: 

\begin{lemma}\label{LQmax}
Let $q>1$ and $h_n=\max_{x\in g_n} v(x)$. Then
\[
{1\/ 2\pi}\sum_{n=1}^{q-1}  h_n|g_n|<Q_0< {1\/ \pi}\sum_{n=1}^{q-1}  h_n|g_n|.
\label{Q<>}
\]
\end{lemma}

\no{\it Proof.} We have
$$
Q_0={1\over\pi}\int_0^\pi v(x)dx<
{1\over\pi}\sum_{n=1}^{q-1}h_n|g_n|,
$$
which is the r.h.s. inequality of (\ref{Q<>}).
To establish the l.h.s. we observe that the concavity of 
$v(x)$ on gaps implies that
$$
  \int_{g_n}v(x)dx>{1\/ 2} h_n|g_n|
$$
and then use (\ref{v0}). $\BBox$

In order to find a bound on the gap width from below, we shall need 
an upper bound on $h_n$:

\begin{lemma}\label{Lh+}
Let $q>1$ and $h_+=\max_n h_n$. Then
\[
0<h_+<\ln\left({c\over A}+|{1\over Aq}\sum_{j=1}^q b_j|
\right)<\ln{2c\over A}.
\]
\end{lemma}

\no{\it Proof}. We begin with an argument similar to that in \ci{Last}.
Let $\l_j$ be the eigenvalues of $L$ (zeros of $D(\l)$).
By the inequality between the algebraic and geometric means
$$
A^q |D(\l)|=\prod_{j=1}^q|\l-\l_j|<(S(\l)/q)^q,
$$
where $S(\l)=\sum_{j=1}^q|\l-\l_j|$. Observe that
for $\l\in[\min\l_n,\max\l_n]$
\[
S(\l)\le \max\left\{\sum_{j=1}^q(\l_j-\min\l_n),
\sum_{j=1}^q(\max\l_n-\l_j)\right\}.\label{S}
\]
Hence, $S(\l)<cq+|\sum_{j=1}^q \l_j|=
cq+|\sum_{j=1}^q b_j|<2cq$, where we used the fact that $H$ is normalised. 
Now note that
$$
h_+={1\over q}{\rm arccosh}\;{|D(\l_+)|\over 2}<
{1\over q}\ln |D(\l_+)|<\ln(S(\l_+)/Aq),
$$
where $\l_+$ is a point of maximum of $|D(\l)|$ in the gaps.
The positivity of $h_+$ follows from Lemma~\ref{Lu}.
$\BBox$

We have
\[
{1\over\pi}\int_0^\pi v(x)dx<
{1\over\pi}\sum_{n=1}^{q-1}h_n|g_n|<{h_+\over\pi}\sum_{n=1}^{q-1}|g_n|,
\]
and, in view of (\ref{v0}), we have the lower bound (\ref{g1})
for the total gap width. The inequality (\ref{g1p}) is proved similarly
using (\ref{v2}).


To prove the r.h.s. of (\ref{g2}), we use the semicircle property
of $v(x)$:
$$
Q_0=\ln{c\over 2A}={1\/ \pi}\int_{0}^{\pi}v(x)dx>
{1\/ \pi}\sum_{n=1}^{q-1}\frac{\pi(|g_n|/2)^2}{2}={1\over 8}
\sum_{n=1}^{q-1}|g_n|^2,
$$
which yields the r.h.s. of (\ref{g2}).

We are yet to prove the l.h.s. of (\ref{g2}).
For this we shall make use of the Dirichlet integral (\ref{Dir1}).

\begin{lemma}\label{4.3}
We have the following estimate:
\[
\sum_{n=1}^{q-1}h_n^2\le\pi^2b_+\ln(c/2A),
\quad
{\it where} \quad
b_+=\max\left\{1,{q h_+\over\pi}\right\}.
\]
\end{lemma}

\no {\it Proof.} 
First, recall the following result from \ci{K3}.
Let a real function
$f(u+iv)$ belong to the Sobolev space $  W^1_2(D(h,\alpha,\beta)),$ 
where the domain
$D(h,\alpha,\beta)=\{ \alpha < u < \beta  , \ 0 < v < h  \}$
for some $\beta > \alpha, h > 0 $.
Let $f$ obey the following conditions :
 
(a) $f(u+i0)=0,  $ if $u\in (\alpha, \beta ),$

(b) $f$ is continuous in the closure $\bar D(h,\alpha,\beta).$ 

\no Then it was shown in \ci{K3} that
\[
\int _0^{h}{|f(+ \alpha +iv) |^2dv\over v}\leq
{\pi \over 2}\max\left\{1,{h\over{\beta-\alpha}}\right\} 
\iint_{D(h,\alpha,\beta)}|\nabla f|^2dudv.\label{5.2}
\]

Now take the function $f(u+iv)={\rm Im}(k-z(k))=v-y(k)$ and note that

(a) the identity $y(u+i0)=0$ yields $f(u+i0)=0$;

(b) as the function $y(k)$ is continuous in $\cK_+^0$, 
$f$ is also continuous there.

\no
Therefore the conditions leading to (\ref{5.2}) are fulfilled
for this $f(k)$ and the domain\newline
$D(h_n,\pi n/q,\pi(n+1)/q)$.
Since $f(\pi n/q+iv)=v$ if $0<v<h_n$, we get
\[
h_n^2=2\int_0^{h_n}vdv \leq \pi
  \max\left\{1,{q h_n\over\pi}\right\}
\iint_{D}|z'(k)-1|^2dudv.
\]
Replacing $\max\left\{1,{q h_n\over\pi}\right\}$ 
by its maximum $b_+$ and summing over $n$, we obtain
\[
\sum_{n=1}^{q-1}h_n^2\le\pi b_+\iint_{\cK_+^0}|z'(k)-1|^2dudv=
\pi b_+\iint_{\cZ_+^0}|k'(z)-1|^2dxdy,
\]
where the last equality follows by change of variables.
In view of (\ref{Dir1}) we finished the proof. $\BBox$

Applying Lemma~\ref{4.3}, we get
$$
\ln{c\over 2A}=
{1\over\pi}\int_0^\pi v(x)dx<
{1\over\pi}\sum_{m=1}^{q-1}h_m|g_m|<
{1\over\pi}\sqrt{\sum_{m=1}^{q-1}h_m^2}\sqrt{\sum_{m=1}^{q-1}|g_m|^2}<
\sqrt{b_+\ln{c\over 2A}}\sqrt{\sum_{m=1}^{q-1}|g_m|^2},
$$
which yields the l.h.s. of (\ref{g2}).

\no\section*{Acknowledgments}
The authors were partly supported by the Sfb288. I.K. is grateful to
S. Jitomirskaya, D. Lehmann, H. Schulz-Baldes, and F. Sobieczky
for useful discussions.

\no

\end{document}